\documentclass{amsart}
\usepackage{amsfonts}
\usepackage{amsmath}
\usepackage{amssymb}
\newtheorem{remark}{Remark}
\def\endproof{{\par \raggedleft \vrule height4pt width 4pt \par}}
\def\proof{{\bf Proof \par}}

\def\poly {{I \!\! P}}

\def\integer {\mbox{\sf Z}\!\!\mbox{\sf Z}}
\def\proj {{\mathcal P}}

\def\Ball{ {\mathcal {B}}}
\def\BigO{ {\mathcal {O}}}
\def\B{\hat \beta}
\def\f{\hat f}

\def\QL{Q^{\lambda}}

\def\sinc{\mathop{\mbox{\rm Sinc}}\nolimits}
\def\supp{\mathop{\mbox{\rm Supp}}\nolimits}

\def\R{\mathbb R}
\def\var{\varepsilon}

\newtheorem{Defi}{Definition}
\newtheorem{Thm}{Theorem}
\newtheorem{Prop}[Thm]{Proposition}
\newtheorem{Cor}{Corollary}[Thm]

\newcommand{\bq}{\begin{equation}}
\newcommand{\eq}{\end{equation}}
\def\bqa{\begin{eqnarray}}
\def\eqa{\end{eqnarray}}
\def\bd{\begin{displaymath}}
\def\ed{\end{displaymath}}

\begin{document}
\title[Spectral methods for the non cut-off Boltzmann equation]{Spectral methods for the non cut-off Boltzmann equation
and numerical grazing collision limit}
\author{Lorenzo Pareschi}
\address{Universit\'a di Ferrara, Dipartimento di
Matematica, Via Machiavelli 35, 44100 Ferrara, Italy.
E-mail: {\tt pareschi@dm.unife.it }} 
\author{Giuseppe Toscani}
\address{Universit\'a di Pavia, Dipartimento di
Matematica, Via Ferrata 1, 27100 Pavia, Italy.
E-mail: {\tt toscani@dimat.unipv.it }}
\author{C\'edric Villani}
\address{Ecole Normale Superieure, Departement de
Matematiques, Rue d'Ulm 45, Paris, France.
E-mail: {\tt villani@dmi.ens.fr }}
%
%
\date{January 30, 2002}
\maketitle

\begin{abstract}
In this paper we study the numerical passage from the spatially homogeneous
Boltzmann equation without cut-off to the Fokker-Planck-Landau equation in the
so-called grazing collision limit. To this aim we derive a
Fourier spectral method for the non cut-off Boltzmann equation
in the spirit of~\cite{PP},\cite{PR}.
We show that the kernel modes that define the spectral method have the correct grazing collision limit providing a consistent spectral
method for the limiting Fokker-Planck-Landau equation. In particular, for small values of the scattering angle, we derive an approximate formula for the kernel
modes of the non cut-off Boltzmann equation which, similarly to the
Fokker-Planck-Landau case, can be computed with a fast algorithm. The uniform spectral accuracy of the method
with respect to the grazing collision parameter is also proved.


\bigskip
\noindent
{\bf Key words}~ Spectral methods, Boltzmann equation, cut-off assumption, Fokker-Planck-Landau equation, grazing collision limit.
\end{abstract}

\bigskip
\maketitle

\section{Introduction}

In this paper we are concerned with numerical resolution methods for collisional equations arising both from kinetic theory of rarefied gases and plasma physics. These equations take the form
\begin{equation}
\partial_tf + v \cdot \nabla_x f = Q(f,f), \qquad t\ge 0, \quad x, v\in \R^3
\label{kinetic}
\end{equation}
where $f(x,v,t)$ is assumed to be a nonnegative function that represents the density of particles in position $x$ at time $t$ with velocity $v$. In (\ref{kinetic}), $Q(f,f)$ is a quadratic integral operator describing the collisions of particles, whose form we shall make explicit below. Here we will mainly focus our attention to the approximation of $Q(f,f)$.

As it is well known in the physics literature, the Fokker-Planck-Landau (FPL) equation is used to describe the binary collisions, occurring in a plasma, between charged particles and can be written as
\begin{equation}
\partial_tf + v \cdot \nabla_x f = Q_L(f,f),
\label{landau}
\end{equation}
\begin{equation}
Q_L(f,f) = \nabla_v \cdot \left ( \int_{\R^3} dv_*\, a(v-v_*)
\left [ f_* \nabla_v f - f (\nabla_{v} f)_* \right ] \right ),
\label{QL}
\end{equation}
where $a$ is the matrix-valued function of the form
\[ a_{ij}(z) = \Psi(|z|) \left ( \delta_{ij} - \frac{z_iz_j}{|z|^2}
\right ),\]
and we use the standard notations $\phi_*=\phi(x,v_*,t)$.
The kernel $\Psi$ depends on the interaction between
particles. A typical choice is, up to a multiplicative constant
\[ \Psi(|z|) = |z|^{\gamma+2}, \qquad -3 \leq \gamma \leq 1. \]

This operator was introduced as an approximation of the Boltzmann collision
operator in the case of {Coulomb interaction} (i.e. $\gamma=-3$)
\cite{landau:Coulomb:36}.
The motivation is that the Boltzmann operator is meaningless in the case
of a Coulomb interaction. In fact, for a Coulomb
interaction, the effect of those collisions that are {\em grazing}, i.e. collisions which result in an infinitesimal angle deflection of the particles trajectories, prevails over the effect of other collisions.

To clarify this idea, let us consider
the Boltzmann equation in the form~\cite{CE},\cite{CIP}
\begin{equation} \label{EB}
 \partial_tf + v \cdot \nabla_x f = Q_B(f,f),
\end{equation}
\begin{equation} \label{BCO}
 Q_B(f,f) = \int_{\R^3}dv_* \int_{S^2}d\sigma\, B(v-v_*,\theta)
\bigl ( f'f'_* - ff_* \bigr ),
\end{equation}
where the deflection angle $\theta \in [0,\pi/2]$, is such that $\cos \theta = (v-v_*) \cdot \sigma /|v-v_*|$. The unit vector $\sigma$ parameterizes the set of all kinematically possible (i.e., those conserving energy and momentum) post-collisional velocities $v'$ and $v'_*$ by
\[ \begin{cases}
\displaystyle
v' = \frac{v+v_*}{2} + \frac{|v-v_*|}{2} \sigma, \\
\\
\displaystyle
v'_*= \frac{v+v_*}{2} - \frac{|v-v_*|}{2} \sigma.
\end{cases} \]
The relative probability of these outcoming velocities depends on
the nature of the interaction between the molecules, and this is
taken into account in the kernel $B$. For molecules which interact
through an $1/r^s$ force law, where $r$ is the distance between
interacting particles, physical arguments show that the natural
choice for the kernel $B$ is
\[
B(v-v_*,\theta) = |v-v_*|^\gamma b(\cos\theta),
\qquad \cos\theta = (\frac{v-v_*}{|v-v_*|}, \sigma),
\]
where $\gamma = (s-5)/(s-1)$, and $\sin \theta b(\cos \theta)$ is
a smooth function except near $\theta=0$, where it presents a
(nonintegrable) singularity of order $(s+1)/(s-1)$. This
singularity corresponds to the grazing collisions and when these
collisions prevail (in a way we shall make precise in the next
section) solutions to the Boltzmann equation converge towards
solution of the FPL equation. In recent years, a noticeable amount
of theoretical
results~\cite{arsen:landau:90},\cite{degond},\cite{desvill:92},
\cite{vill:new:98},\cite{vill:lm:98} concerned with the limiting
process involved, called the {grazing collisions asymptotics}, has
been obtained.

On the other hand, the usual numerical resolution methods for the
Boltzmann equation can not be used in this context, exactly in
reason of the singularity in the kernel $B$. Classically, the
standard strategy for numerical applications is to {\em cut-off}
the small angle collisions so that $B$ becomes integrable. In
fact, under this simplification, it is possible to split
$Q_B(f,f)$ into a gain term $Q^+_B(f,f)$ and a loss term
$Q^-_B(f,f,)$ according to
\[
Q_B(f,f) = Q^+_B(f,f) - Q^-_B(f,f),
\]
which makes the problem treatable both by direct simulation Monte
Carlo methods~\cite{bird},\cite{ilne},\cite{nambu} or fully
deterministic methods~\cite{buet},\cite{rosc}. In particular, as
the cut-off parameter becomes smaller, the computational cost of
these methods increases dramatically so that they are unable to
treat in practice the non cut-off problem. Recently a study of the
behavior of Monte Carlo methods for a simplified Boltzmann model
close to the non cut-off case has been presented in
\cite{desvill:99}.

To the best of our knowledge, the method we have studied in the
present paper is the first method which is able to treat
numerically the Boltzmann collision operator in the non cut-off
case.

In very recent years, there has been a new approach, based on a
Fourier spectral method, to the numerical approximation of the
Boltzmann and the FPL operators~\cite{PP},\cite{PR},\cite{PRT}.
The advantage of this method is that the whole structure of the
collisional operator (in practice all the informations
characterizing the kinetic equation) is contained into a series of
kernel modes. In addition, in the FPL case, thanks to the
structure of these kernel modes, the resulting algorithm can be
evaluated with a strong reduction in the computational cost. Other
numerical methods, that avoid the quadratic complexity of the FPL
operator, have been introduced recently
\cite{bcdl},\cite{dl},\cite{Lemou}.

Here, we consider this spectral method in the case of the Boltzmann collisional operator (\ref{BCO}), without the cut-off assumption, by showing that the corresponding kernel modes are well-defined for any reasonable kernel with singularity. The main result here is that these kernel modes converge (in the grazing limit) to the corresponding kernel modes of the FPL equation and that the spectral method is uniformly accurate with respect to the grazing collision parameter.
In addition, we derive an approximate representation for the kernel modes of the non cut-off Boltzmann equation which is valid for small values of the deflection angle and that, similarly to the FPL case, can be computed through a fast algorithm. Numerical results based on this representation are actually under development and will be presented elsewhere.

Summarizing the method has the following properties:
\begin{itemize}
\item Well-defined for the non cut-off Boltzmann equation.
\item Correct numerical grazing collision limit.
\item Uniform spectral accuracy with respect to the grazing collision parameter.
\item Strong reduction of the computational cost close to the grazing limit.
\end{itemize}

For the sake of completeness, we mention here that the idea of using the numerical asymptotic limit as a guide
for the development of numerical methods is related to some recent works in the context of the fluid-limit of the Boltzmann equation (see~\cite{CJR},\cite{GPT} and the references therein).

 We briefly discuss the organization of the paper. Section 2 is devoted to recall some mathematical and physical properties of the equations and to present the precise formulation of the grazing collision limit. Section 3 deals with the derivation of the spectral projection of the non cut-off Boltzmann equation. The grazing collision limit of the kernel modes is contained in Section 4.
Finally, Section 5 is devoted to a proof of the uniform spectral accuracy of the method. Some final considerations are contained in Section 6.

\section{Mathematical and physical background}

In this section, we briefly recall the mathematical and physical background both of the  FPL and Boltzmann equations.
In the rest of the paper we will assume a kernel $B$ with a non integrable singularity, as described in the introduction, and we will consider the space homogeneous equations. It is well-known, in fact, that by the standard splitting algorithm we may restrict to the space homogeneous case.

From the physical point of view, the most important quantity
associated to the angular cross-section $b(\cos\theta)$ is the so-called
(angular) cross-section for momentum transfer, which is (up to a
multiplicative constant)
\[ \Lambda = \int_{S^{N-1}} b(k\cdot\sigma) (1-k\cdot\sigma)\,d\sigma  =
2\pi \int_0^{\pi/2} b(\cos \theta)
(1-\cos\theta)\, \sin\theta\, d\theta, \]
or more generally the function
\begin{equation}
\label{mt}
A(q)=2\pi \int_0^{\pi/2} B(q,\theta) (1-\cos \theta)\,\sin\theta\, d\theta.
\end{equation}
It is shown in~\cite{vill:new:98} that for $-3 \leq \gamma \leq 1$,
the finiteness of $\Lambda$ is exactly what is needed to develop a
mathematical existence theory.

If $s=2$ (Coulomb interaction), then $\Lambda$ is a divergent integral
as $\theta \to 0$, which prevents to define $Q_B(f,f)$ in any
reasonable sense (see~\cite[Part I, Appendix A]{vill:thesis}).
Using heuristic arguments, Landau showed that the principal
part of~\eqref{BCO} under the truncation $(\theta \geq \varepsilon >0)$
is proportional to $Q_L(f,f)$ (see~\cite{vill:div:97} for a rigorous
variant of Landau's argument).

In view
of these remarks, the following setting will appear
to be natural from both the mathematical and physical points of view.
We use the standard notations
\[ L^1_s = \left \{ f \in L^1(\R^3), \quad \int_{\R^3} |f(v)|(1+|v|^s)\,dv <
\infty \right \}, \]
\[ L \log L = \left \{ f \in L^1(\R^3), \quad \int_{\R^3} |f(v)|
\log (1+|f(v)|)\, dv < \infty \right \}. \]
Moreover, if $f$ is a nonnegative function, we shall
call $\int f\,dv$ its mass, $\int fv\,dv$ its momentum,
$\int f(|v|^2/2)\,dv$ its energy, and $\int f \log f\,dv$ its entropy.

\begin{Defi}
We shall say that the family of kernels $(B_\var)_{\var>0}$
given by
\begin{equation} \label{family}
B_\var(q,\theta) = |q|^\gamma b_\var(\cos\theta)
\end{equation}
($-3\leq \gamma \leq 1$) concentrates on grazing collisions if
\[ a)\qquad \Lambda_\var \equiv
2\pi \int_0^{\frac{\pi}{2}} b_\var(\cos\theta)
(1-\cos\theta)\, \sin\theta\, d\theta
\longrightarrow   \Lambda_0 < \infty, \]
\[ b) \qquad \forall \theta_1>0, \qquad
b_\var(\cos\theta) \longrightarrow 0
\quad\text{uniformly on $(\theta \geq \theta_1)$}. \]
\end{Defi}

In fact, we could allow much more general assumptions on
$B_\var$. For instance,
\begin{description}
\item[(i)]
\[
\int_0^{\pi/2} B_\var(q,\theta) \sin^3 \theta\,d\theta
 \leq C \left ( \frac{1}{|q|^3} + |q| \right )
\]
\item[(ii)]
\[
2\pi
|q|^2 \int_0^{\pi/2} B_\var(q,\theta) (1-\cos\theta)\,
\sin\theta\,d\theta \rightharpoonup A_0(|q|)
\]
in $W^{1,\infty}_{\text{loc}}(\R^3 \setminus \{0\})$,
\item[(iii)]
\[
\forall R>0, \qquad
\int_{|q|\leq R} |q|^2 \int_0^{\pi/2} B_\var(q,\theta) \sin^4 \theta\,
d\theta\,dq \longrightarrow 0.
\]
\end{description}

And then, the function $A_0(|q|)$ would replace the function $\Lambda_0
|q|^{2+\gamma}$ in all that follows. The main results of~\cite{vill:new:98} can then be stated
as follows,

\begin{Thm} Let $f_0$ be an initial nonnegative datum with finite mass,
energy and entropy. If $\gamma>0$, assume in addition that $f_0 \in
L^1_{2+\delta}$ for some $\delta>0$. Then, for all $\var>0$ there
exists a solution $f^\var$ of the Boltzmann equation with kernel
$B_\var$, such that $f^\var(0,\cdot)=f_0$. Moreover, as $\var \to 0$,
up to extraction of a subsequence, for all $T>0$ $f^\var$ converges weakly in
$L^1([0,T];\R^3)$ to a solution $f$ of the FPL equation with
kernel $\Psi(|z|)=\frac18 \Lambda|z|^{\gamma+2}$.
\end{Thm}

\begin{remark}
\rm The exact value of the cross-section given in the previous theorem, comes out in a clear way from the analysis of the numerical grazing collision limit  we will develop in Section 4.
\end{remark}

We do not make precise here the meaning of weak solutions, which is not
the usual one for $\gamma < -2$ (see~\cite{vill:new:98} for details).

From now on, we use the notation
\[ \zeta(\theta) = b(\cos\theta)\, \sin\theta, \]
in which the Jacobian arising from the spherical variables is
taken into account, and we recall that this function has typically
a singularity at $\theta=0$ like $\theta^{-(1+\nu)}$,
$\nu = 2/(s+1)$ for inverse $s$-forces.
Typical families $(\zeta_\var)$ are given by
\[ \begin{cases}
\displaystyle
\zeta_\var(\theta) = \frac{1}{\ln \var^{-1}} \zeta(\theta)1_{\theta
\geq \var} \quad \text{for $s=2$}, \\
\\
\displaystyle
\zeta_\var(\theta) \sin^2 \frac{\theta}{2} =
\frac{1}{\theta}\sin^2 \frac{\theta}{2\var}
\zeta(\frac{\theta}{\var}) \quad \text{for $s>2$}.
\end{cases} \]

We emphasize that from the physical point of view, the case
$\gamma=-3$ is the most interesting (see also~\cite{spohn:large:91},
in which a somewhat different viewpoint on the FPL equation is
given).

Now, let us give the state of the art concerning the mathematical
theory of the spatially homogeneous FPL equation. The first works
on the subject were done by Arsen'ev and
Buryak~\cite{arsen:landau:90}. Their results were considerably
improved and extended in the recent work~\cite{desvvill:II:97}, from
which we extract the following theorem, which gives a quite
satisfactory picture of the case $\gamma>0$~

\
\begin{Thm} Let $\gamma>0$, and let $f_0$ be a
nonnegative initial datum, $f_0 \in L^1_{2+\delta}$ for some
$\delta>0$. Then, there exists a solution $f(t,v)$ of the FPL
equation, continuous in time,
with $f(0,\cdot) = f_0$, satisfying the
conservation of mass, energy, and the decrease of entropy. and
which is very smooth for all
positive time. More precisely,
for all $t \geq t_0>0$, and $k,s \geq 0$,
one has $f \in H^k \cap L^1_s$, where $H^k$
denotes the standard Sobolev space of order $k$, and the norm of $f$
in these spaces depends only on $f_0$ and $t_0$.
In addition, if $\int |f_0(v)|^2 (1+|v|^s) <\infty$,
with $s> 15+5\gamma$, then $f$ is unique (in the class of weak
solutions which may let the energy decrease) and for all time $t>0$, there
exist positive constants $A_t$ and $K_t$ such that $f(t,v) \geq A_t
e^{-K_t|v|^2}$.
\end{Thm}

In the case $\gamma=0$, it is sufficient to require that $f \in L^1_2$
to get both existence and uniqueness of solutions~\cite{vill:lm:98}.

As regards other qualitative features of solutions, a detailed study
of the asymptotic behavior as $t \to \infty$ was done
in~\cite{desvvill:II:97}, which led to the following results

\
\begin{Thm} Assume $\gamma>0$, and
let $f$ be a (smooth) solution of the FPL equation with
kernel $|z|^{\gamma+2}$. Without loss of generality, assume
 that $f$ has mass~1,
momentum~0 and energy $3/2$. Then,
\[ \|f(t,\cdot) -M \|_{L^1} \leq C(f_0)\, t^{-1/\gamma}, \]
where $C(f_0)$ is explicit and depends only on the entropy of $f_0$,
and
\[ M(v) = \frac{e^{-\frac{|v|^2}{2}}}{(2\pi)^{3/2}}. \]
Moreover, $f-M$ goes to 0 in all (weighted) Sobolev norms.
\end{Thm}

Again, in the case $\gamma=0$ the results are simpler, in that the
trend to equilibrium is exponentially fast, with a bound of the form
$C_\var e^{-(2-\var)t}$, for all $\var >0$.

As a consequence of these results, it was shown
in~\cite{desvvill:II:97} that if $\gamma \geq 0$, and if
the initial datum satisfies the conditions for uniqueness of the
solution, then this solution is stable, {\em globally in time},
say in $L^1$ norm, with respect to small perturbations of both the
initial datum (in a suitable weighted $L^2$ space) and the
kernel, in the sense that one allows perturbations of the form
$\Psi(|z|) = |z|^{\gamma+2}(1+\eta(|z|))$, with $\eta$ small in $C^2$
norm. From the practical viewpoint, this means that small errors in
computations are allowed.

We conclude this presentation by noting that in the case $\gamma <0$,
the mathematical theory is still at an early stage (only existence of
global weak solutions is known, and the study of
trend to equilibrium is still in progress).

\section{Spectral projection of the non cut-off Boltzmann equation}
In this section we derive the Fourier spectral method for the non cut-off
Boltzmann equation following~\cite{PP},\cite{PR}.

A simple change of variables permits to write the Boltzmann collision
operator $Q_B(f,f)$ in the form
\begin{equation}
Q_B(f,f) = \int_{\R^3}dq \int_{S^2} d\sigma
B(q, \theta)[f(v+q^+) f(v+q^-) - f(v) f(v+q)],
\label{eq:QBG}
\end{equation}
where $q = v_\ast-v$ and the vectors $q^+$ and
$q^-$ that parameterize the  post-collisional velocities are given by
\begin{equation}
q^{+} = \frac12(q+\vert q\vert \sigma), \quad
q^{-} = \frac12(q-\vert q\vert \sigma).
\label{eq:VV2}
\end{equation}
We point out that the possibility to integrate the collision operator
over the relative velocity is essential in the derivation of the method.
Other kinetic equations with a similar structure may be treated in the same way.

First we need to reduce the problem to a bounded domain in $v$. To this aim
we observe that if a distribution function $f$ has compact support,
$\supp (f(v)) \subset \Ball(0,R)$, where $\Ball(0,R)$ is the ball of
radius $R$ centered in the origin, then by conservation of energy and
momentum the collision operator  satisfies the following \cite{PP}

\begin{Prop}
Let $\supp (f(v)) \subset \Ball(0,R)$ then
\[
Q_B(f,f)(v) = \int_{\Ball(0,2R)} dq\,\int_{S^2} d\sigma
B(\vert q \vert, \theta)[f (v+q^+) f (v+q^-) - f(v) f(v+q)]\,,
\]
with $v+q^+,v+q^-,v+q \in \Ball(0, ( 2 + \sqrt{2} )R)$.
\label{pr:1}
\end{Prop}
Thus to develop a spectral method we
will consider the distribution function $f(v)$ restricted on the cube
$[-T,T]^3$ with $T \geq (2+{\sqrt 2})R$, assuming $f(v)=0$ on $[-T,T]^3
\setminus \Ball(0,R)$,
and extend it by periodicity to a periodic function on $[-T,T]^3$.

As observed in \cite{PR} using the periodicity of the function it is
enough to take $T \geq (3+{\sqrt 2})R/2$ to prevent intersections of
the regions where $f$ is different from zero.

To simplify the notation let us take $T=\pi$ and hence $R=\lambda\pi$
with $\lambda = 2/(3+\sqrt{2})$. Hereafter, we use just one index to
denote the three-dimensional sums with respect to the vector
$k=(k_1,k_2,k_3)\in \integer^3$.

The approximate function $f_N$ is represented as the truncated Fourier
series
\begin{equation}
f_N(v) = \sum_{k=-N}^N \f_k e^{i k \cdot v},
\label{eq:FU}
\end{equation}
\begin{equation}
\f_k = \frac{1}{(2\pi)^3}\int_{[-\pi,\pi]^3} f(v)
e^{-i k \cdot v }\,dv.
\label{eq:FC}
\end{equation}
The Fourier-Galerkin method~\cite{CHQZ},\cite{GO} is obtained by considering the
projection of the homogeneous Boltzmann equation on the space of trigonometric
polynomials of degree $\leq N$.

Hence, we have
\begin{equation}
\int_{[-\pi,\pi]^3}
\left[{\partial_t f_N} - Q_B(f_N,f_N)
\right]
e^{-i k \cdot v}\,dv = 0, \quad k=-N,\ldots,N.
\label{eq:VAR}
\end{equation}
By substituting expression (\ref{eq:FU}) in (\ref{eq:VAR}) we get
a set of ordinary differential equations for the Fourier coefficients
\begin{equation}
{\partial_t \f_k} = \sum_{\underset{l+m=k}{l,m=-N}}^N \f_l\,\f_m
\B(l,m),\quad k=-N,\ldots,N,
\label{eq:CF1}
\end{equation}
where the {\em Boltzmann kernel modes} $\B(l,m)$ are given by
\begin{equation}
\B(l,m) = \int_{\Ball(0,2\lambda\pi)} dq\, \int_{S^2} d\sigma\,
B(q, \theta) \left[e^{il\cdot q^+ + im\cdot q^-} - e^{i q \cdot m }\right].
\label{eq:KM}
\end{equation}

Note that, as a consequence of the properties of
trigonometric polynomials, (\ref{eq:KM}) is a scalar quantity
completely independent on the argument $v$,
depending just on the particular kernel structure. This property is
strictly related to the use of a Fourier spectral method. Other
spectral methods may be developed, however they do not
lead to this fundamental simplification.

In practice all the
informations characterizing the kinetic equation are now contained in the
kernel modes. Clearly, these quantities can be computed in advance and then
stored in a multidimensional matrix.

Finally, the integral in (\ref{eq:KM}) can be written in the form
\[
\B(l,m) =
\int_0^{2\lambda\pi} \int_{S^2}
\int_{S^2} B(\rho, \sigma\cdot\omega) e^{i\rho(|l+m|\omega\cdot \eta)}
\left[e^{-i\rho(|l-m|\sigma\cdot\mu)}
      - e^{-i\rho(|l-m|\omega\cdot\mu)}\right]
\rho d\rho\,d\omega\,d\sigma,
\]
where $q=\rho\omega$, $(l+m)=|l+m|\eta$ and $(l-m)=|l-m|\mu$.

It is possible to show that the previous expression is invariant under
any unitary transform $R$ (rotation or symmetry) of the unit
vectors $\eta$ and $\mu$ \cite{PP}. As a consequence the kernel modes
are functions only of $|l-m|$, $|l+m|$ and of the angle between
$\eta$ and $\mu$. In particular we have
\begin{equation}
\B(l,m)=F(|l+m|,|l-m|,\eta\cdot\mu)=F(|l-m|,|l+m|,\eta\cdot\mu).
\label{eq:id1}
\end{equation}
Since from (\ref{eq:KM}) it follows easily
$\B(l,m)=\overline{\B(-l,-m)}$, from (\ref{eq:id1}) we get also
$\B(l,m)=\B(-l,-m)$ and hence the kernel modes are real functions.

Obviously, these properties are useful to reduce the storage requirements
of the method. We refer to \cite{PR} for a more detailed
discussion on this topic.

\begin{remark}
\rm
\begin{itemize}
\item[]
\item[i)]
The evaluation of (\ref{eq:CF1}) requires exactly $\BigO (N^6)$
operations. We emphasize that the usual cost for a method based on $N^3$
parameters for $f$ in the velocity space is $\BigO (n_a N^6)$ where $n_a$
is the numbers of angles discretizations. Thus the straightforward
evaluation of (\ref{eq:CF1}) it is less expensive than a usual
discrete-velocity algorithm.

\item[ii)] In the VHS case, $B(q,\theta)=C_\alpha |q|^{\alpha}$, the dependence on the
scattering angle disappears and
(\ref{eq:KM}) reduces to a one-dimensional integral \cite{PR}
\[
\B(l,m) = C_{\lambda,\alpha} \int_0^{\lambda\pi} r^{2+\alpha}\left[
\sinc(|l+m|r)
\sinc(|l-m|r)-\sinc(2|m|r)\right]\,dr,
\]
where $C_{\lambda,\alpha}=(4\pi)^2(8)^{\alpha}C_\alpha$.
In addition for integer values of $\alpha < -3$ (like in the case
of Maxwell molecules or hard spheres) the previous integral
can be computed explicitly \cite{PR} given rise
to very simple and easy implementable approximation formulas.
\end{itemize}
\end{remark}

\section{The grazing collision limit of the kernel modes}
From the previous section it is clear that in the grazing collision limit
it would be desirable that the previous spectral projection gives a
consistent spectral method for the limit FPL equation.

To this aim first we compute the grazing collision limit of the Boltzmann
kernel modes and then we verify directly on the FPL equation that the
asymptotic limit of the Boltzmann kernel modes coincides with the
kernel modes of the FPL equation.

\subsection{Computation of the grazing collision limit}
We shall compute the limit of the Boltzmann kernel modes, as the parameter
$\epsilon$ defining the asymptotics of grazing collision goes to $0$.

Thus we start from (\ref{eq:KM}) written in the form
\begin{equation}
\B_\var(l,m) = \int_{\Ball(0,2\lambda\pi)} dq\, e ^{iq\cdot m} \int_{S^2} d\sigma\,
B_\var(q, \theta) \left[e^{i\frac{(l-m)}{2}\cdot (q-|q|\sigma)} - 1\right],
\label{eq:KM2}
\end{equation}
where $B_\var(q,\theta)$ is given by (\ref{family}).

Let us introduce a system of spherical coordinates around the vector $q$. Hence
$(q/|q|,j,h)$ will be a fixed basis, and the spherical angle $\phi$ will be
defined by
\[
|q|\sigma = q\cos\theta + j |q|\sin\theta\sin\phi + h|q|\sin\theta\cos\phi.
\]
We get
\[
\B_\var(l,m) = \int_{\Ball(0,2\lambda\pi)} dq\, e^{iq\cdot m} \int_{0}^{\pi}
d\theta\sin\theta\, B_\var(q, \theta)
\int_{0}^{2\pi} d\phi\,[e^{F_{l-m}(q,\theta,\phi)}-1],
\]
with
\begin{equation}
F_{l-m}(q,\theta,\phi)=i\frac{(l-m)}{2}\cdot [q (1-\cos\theta)+
j|q|\sin\theta\sin\phi + h|q|\sin\theta\cos\phi].
\end{equation}
To start with, we use the Taylor expansion
\begin{equation}
e^{F_{l-m}(q,\theta,\phi)}-1 = F_{l-m}(q,\theta,\phi)+\frac12 [F_{l-m}(q,\theta,\phi)]^2 + \BigO(|q|^3|l-m|^3\theta^3).
\label{taylor}
\end{equation}
By the assumptions in Definition~1 (or point $(iii)$
after Definition 1), the remainder can be
neglected since as $\var \to 0$ we have
\[
\int_{\Ball(0,2\lambda\pi)} dq\, e^{iq\cdot m} |q|^3 \int_{0}^{\pi}
d\theta\sin\theta\, B_\var(q, \theta) \theta^3 \to 0.
\]
Now we consider the contribution of the first order term in (\ref{taylor}).
By symmetry we have that
\[
i|q|\int_0^{2\pi}d\phi\,\frac{(l-m)}{2}\cdot [j\sin\theta\sin\phi +
h\sin\theta\cos\phi] = 0.
\]
Thus the first order contribution reduces to
\begin{equation}
\pi i \int_{\Ball(0,2\lambda\pi)} dq\, e^{iq\cdot m}
(l-m)\cdot q
\int_{0}^{\pi}
d\theta B_\var(q, \theta)\sin\theta (1-\cos\theta).
\label{first}
\end{equation}
Next we consider the second order term in (\ref{taylor}). Again by simple
symmetry considerations the terms in $\sin^2\theta\sin\phi\cos\phi$, $\sin\theta\sin\phi$ and $\sin\theta\cos\phi$ disappear. Moreover also
\[
-\frac{\pi}{4} \int_{\Ball(0,2\lambda\pi)} dq\, e^{iq\cdot m}
[(l-m)\cdot q]^2
\int_{0}^{\pi}
d\theta B_\var(q, \theta)\sin\theta (1-\cos\theta)^2,
\]
is small of order $\var^2$ since $(1-\cos\theta)^2$ is of order $\theta^4$.

For the remaining terms we have
\[
-\frac{1}{8}|q|^2\sin^2\theta \int_0^{2\pi}d\phi\,
\{[(l-m)\cdot j]^2\sin^2\phi +
[(l-m)\cdot h]^2\cos^2\phi\}
\]
\[
= -\frac{\pi}{8}|q|^2\sin^2\theta \{[(l-m)\cdot j]^2 + [(l-m)\cdot h]^2\}
\]
\[
= -\frac{\pi}{8}|q|^2\sin^2\theta [(l-m)^\perp]^2,
\]
where $(l-m)^\perp = (l-m) - ((l-m) \cdot q /|q|) q / |q|$.

Thus the second order contribution is given by
\begin{equation}
-\frac{\pi}{8} \int_{\Ball(0,2\lambda\pi)} dq\, e^{iq\cdot m}
|q|^2 [(l-m)^\perp]^2
\int_{0}^{\pi}
d\theta B_\var(q, \theta)\sin^3\theta.
\label{second}
\end{equation}

Finally, adding (\ref{first}) and (\ref{second}) we obtain the approximation for
small $\var$
\begin{equation}\label{approx}
\B_\var(l,m) = 2\pi\int_{\Ball(0,2\lambda\pi)} dq\, e^{iq\cdot m}
\int_{0}^{\pi}
d\theta\sin\theta \sin^2\frac{\theta}{2} B_\var(q, \theta)\,
G_{l-m}(q,\theta) +\BigO(\var),
\end{equation}
with
\begin{equation}
G_{l-m}(q,\theta) =
i q \cdot (l-m)-\frac14\cos^2\frac{\theta}{2}|q|^2 [(l-m)^\perp]^2.
\label{G}
\end{equation}
Since by assumptions in Definition 1 (or point~(ii) after
Definition~1), we have as $\var \to 0$
\[
2\pi \sin\theta\, (1-\cos\theta) B_\var(q,\theta) \to
\Lambda_0 |q|^\gamma \delta_{\theta=0},
\]
in the limit we obtain the following expression for the
{\em grazing kernel modes}
\begin{equation}
\B_0(l,m) = \frac{\Lambda_0}{2}\int_{\Ball(0,2\lambda\pi)} dq\, e^{iq\cdot m} |q|^\gamma
\{i q \cdot (l-m)-\frac14 |q|^2 [(l-m)^\perp]^2 \}.
\label{grazing}
\end{equation}
Now in (\ref{grazing}) we have a balance between an
advection term, the one with $i q \cdot (l-m)$ and a diffusive term, the one
with $|q|^2[(l-m)^\perp]^2$.

Taking into account the symmetries of the kernel modes from (\ref{eq:id1})
we have $\B_\var(l,m)=\B_\var(-l,m)$ that in the limit $\var \to 0$ gives
the equivalent representation
\begin{equation}
\B_0(l,m) = - \frac{\Lambda_0}{2}
\int_{\Ball(0,2\lambda\pi)} dq\, e^{iq\cdot m}
|q|^\gamma
\{i q \cdot k + \frac14 |q|^2 [k^\perp]^2 \},
\label{grazing2}
\end{equation}
where $k=l+m$.

\begin{remark}
\rm
From the previous results some remarks are in order.
\begin{itemize}
\item[i)] The Boltzmann kernel modes
(\ref{eq:KM2}) are bounded whenever the family of kernels $B_\epsilon$ satisfy the conditions of Definition $1$. In particular, taking $\epsilon =1$, we proved that the Boltzmann kernel modes
are bounded for any kernel $B$ with finite angular cross-section for momentum transfer.

\item[ii)] The final expression (\ref{grazing2}) is much simpler than the starting
expression (\ref{eq:KM2}). It is possible to show \cite{PRT} that this simplification allows to compute the resulting final algorithm in only $\BigO(N^3\log N)$ instead of $\BigO(N^6)$ operations (see section 4.3).

\item[iii)]
Similarly, also the approximation defined by (\ref{approx}-\ref{G}),
that can be used to study numerically the behavior of the non cut-off Boltzmann equation when
the collisions become grazing, can be computed with only $\BigO(N^3\log N)$
operations.
\end{itemize}
\end{remark}

\subsection{Direct derivation of the FPL kernel modes}
We show here that the grazing kernel modes computed in the previous
paragraph coincides with the FPL kernel modes that we can derive by direct
application of the spectral method to the FPL equation.

To start with, for any test function $\varphi$ we have the identity
\[
\int_{\R^3} Q_L(f,f)\varphi(v)\,dv
= -\frac12 \int_{\R^3}dv \int_{\R^3} dv_*\, a(v -v_*)\cdot (\nabla -\nabla_*)
f f_* [ \nabla \varphi - \nabla \varphi_* ]
\]
\[
= \int_{\R^3}dv \int_{\R^3} dv_*\, f f_*\left\{b(v -v_*) \cdot [ \nabla\varphi - \nabla\varphi_* ] +
a_{ij}(v -v_*)\frac12 [\partial_{ij}\varphi + \partial_{ij}\varphi_* ] \right\},
\]
where
\[
b_i(v -v_*)=\sum_j \partial_i a_{ij}(v -v_*).
\]
By symmetry we get
\[
\int_{\R^3} Q_L(f,f)\varphi(v)\,dv =
\int_{\R^3}dv \int_{\R^3} dv_*\, f f_*\left[2b(v -v_*) \cdot \nabla\varphi  +
a_{ij}(v -v_*)\partial_{ij}\varphi \right].
\]
From
\[
a_{ij}(z)= \Psi(|z|) \left(\delta_{ij}-\frac{z_i z_j}{|z|^2}\right),
\]
we have $b_i(z)=-2 z \frac{\Psi(|z|)}{|z|^2}$
(in dimension $n\neq 3$ the factor $2$ has to be replaced by $(n-1)$).

Thus the kernel modes are obtained from
\[
\int_{[-\pi,\pi]^3} Q_L(e^{il\cdot v},e^{im\cdot v})e^{-ik\cdot v}\,dv =
\int_{\Ball(0,2\lambda\pi)} dq\, e^{il\cdot v+im\cdot v_*}
\frac{\Psi(|q|)}{|q|^2}
\{-4 i q \cdot k-|q|^2 [k^\perp]^2 \}e^{-ik\cdot v},
\]
with $k^\perp = k - k \cdot (q/|q|)(q/|q|)$.

Finally we get the expression of the {\em FPL kernel modes}
\begin{equation}
\B_L(l,m) = - 4  \int_{\Ball(0,2\lambda\pi)} dq\, e^{iq\cdot m}
\frac{\Psi(|q|)}{|q|^2} \{i q \cdot k + \frac14 |q|^2 [k^\perp]^2\},
\label{clandau2}
\end{equation}
which is the same  as (\ref{grazing}) for
\begin{equation}
\Psi(|q|) = \frac{\Lambda_0}{8} |q|^{\gamma+2}.
\label{eq:psi}
\end{equation}

Finally, we have proved the following

\begin{Prop} Let $\B_\var(l,m)$ be the kernel modes of the Boltzmann
equation with kernel
$B_\var(q,\theta)=|q|^\gamma \zeta_\var(\theta)$,
where $2\pi \zeta_\var(\theta)(1-\cos\theta)$ is an approximation of
$\Lambda_0 \delta_{\theta=0}$ in the sense of Definition 1. Then we have
\[
\lim_{\var \to 0} \B_\var(l,m) = \B_L(l,m),
\]
where $\B_L(l,m)$ are the kernel modes of the FPL equation with kernel given by (\ref{eq:psi}).
\end{Prop}

\begin{remark}
\rm
\begin{itemize}
\item[]
\item[i)] In arbitrary dimension $n \neq 3$ and with a cross-section like
$\Phi(|q|)$ instead of $|q|^\gamma$, the formula is
\[ \Psi(|q|)=\frac{\Lambda_0}{4(n-1)} \Phi(|q|)|q|^2, \]
as shown under extreme generality in~\cite{av:99}.
We note that the function $\zeta(\theta)$ should of course
now be equal to $b(\cos\theta) \,\sin^{n-2}\theta$, and
that the constant in front of the definition of $\Lambda$,
$\Lambda_0$ should be $|S^{n-2}|$ instead of $2\pi$.

\item[ii)]
It is possible to obtain several representations
of the kernel modes accordingly to the several possible ways in which the FPL equation can be written. For example in \cite{PRT} the FPL kernel modes
has been derived using the same technique of Section 3 in the case of the FPL equation. This originates the expression
\[
\B_L(l,m) =
\int_{\Ball(0,2\lambda\pi)} dq\,
\Psi(|q|) \left\{[l^\perp]^2-[m^\perp]^2\right\}\,
e^{iq\cdot m},
\]
where $l^\perp$ and $m^\perp$ are defined as before.

The equivalence of these representations can be proved by an indirect argument using the recent results in \cite{vill:thesis}.

\item[iii)] Using the symmetry properties of the Boltzmann kernel modes (\ref{eq:id1}) which are conserved in the grazing limit procedure,  we obtain an almost complete set of symmetries for the FPL kernel modes and that also the FPL kernel modes are real functions.
\end{itemize}
\end{remark}

\subsection{Fast algorithms}
We will show here that thanks to the particular structure of equations (\ref{grazing2})
and (\ref{approx}-\ref{G}), that define the FPL kernel modes and an approximation of the
Boltzmann kernel modes in the grazing limit, the resulting spectral method can be evaluated
with $\BigO(N^3\log N)$ operations instead of $\BigO(N^6)$.

We describe the method only for the FPL kernel modes, the extension to expressions of the form
(\ref{approx}-\ref{G}) will then follow easily.

First we rewrite the resulting spectral scheme in the form,
$k=-N,\ldots,N$
\[
{\partial_t \f_k}= \sum_{m=-N}^N \f_{k-m}\,\f_m \B_L(k-m,m),
\]
where we assume that the Fourier coefficients are extended to zero for $|k_j|>N$, $j=1,2,3$.

Next we observe that the term $\B_L(l,m)$ splits as
\[
\B_L(l,m) = - \sum_{j=1}^3 k_j \int_{\Ball(0,2\lambda\pi)} i\,q_j \Psi(q) e^{iq\cdot m}dq
- k^2 \int_{\Ball(0,2\lambda\pi)}|q|^2 \Psi(q) e^{iq\cdot m}dq
\]
\begin{equation}
+ \sum_{j,h=1}^{3} k_j\, k_h  \int_{\Ball(0,2\lambda\pi)} \Psi(q) q_j\, q_h e^{iq\cdot m}dq
\label{eq:split}
\end{equation}
\[
:= \sum_{j=1}^3 k_j F_j(m) + k^2 G(m) +  \sum_{j,h=1}^{3} k_j\, k_h I_{jh}(m),
\]
where $I_{jh}$ is a symmetric matrix.

Thus the resulting scheme requires only the evaluation of convolution sums since it can be written as
\[
{\partial_t \f_k} =
\sum_{j}^{3} k_j \sum_{m=-N}^N \f_{k-m}\,\f_m F_{j}(m) +
k^2 \sum_{m=-N}^N \f_{k-m}\,\f_m G(m)
\]
\[
+ \sum_{j,h=1}^{3} k_j\,k_h \sum_{m=-N}^N \f_{k-m}\,\f_m I_{jh}(m).
\]
It is well-known that transform methods enable
to evaluate the previous expression in
only $\BigO(N^3\log N)$ operations~\cite{CHQZ},\cite{GO}. The computation of the
terms $F_j(m)$, $G(m)$ and $I_{jh}(m)$ can be done following the same
lines presented in \cite{PRT}.

\section{Uniform spectral accuracy of the method}
We explain in this section why the method in this paper
is spectrally accurate, uniformly in the approximation
of grazing collisions.

Let us first introduce some notations. For any
$t \geq 0$, $f_N(v,t)$ is a trigonometric polynomial of degree $N$ in
$v$, i.e. $f_N(t) \in \poly^N$ where
\[
\poly^N = span\left\{e^{ik\cdot v}\,|\, -N \leq k_j \leq N,\, j=1,2,3
\right\}.
\]
Moreover, let $\proj_N : L^2([-\pi,\pi]^3) \rightarrow \poly^N$ be the
orthogonal projection upon $\poly^N$ in the inner product of
$L^2([-\pi,\pi]^3)$  (see (\ref{eq:VAR}))
\[
<f-\proj_N f,\phi>=0,\qquad \forall\,\, \phi\,\in\,\poly^N.
\]
We denote the $L^2([-\pi,\pi]^3)$-norm by $||f||_2 = (< f, f>)^{1/2}.$
In addition
we define $H^r_p ([-\pi,\pi]^3)$ where $r \geq 0$ is an integer to be the subspace of the Sobolev space
$H^r([-\pi,\pi]^3)$, which consists of periodic functions \cite{CHQZ}.

Using these notations the method defined by equation (\ref{eq:CF1}) can be written in equivalent form as
\begin{equation}
\partial_t f_N  = \QL_N(f_N,f_N)
\label{eq:SM}
\end{equation}
with the initial condition $f_N(v,t=0) = f_{0,N}(v)$,
where we have used $\QL(f,f)$ instead of
$\QL_B(f,f)$ to denote the Boltzmann collision operator with relative velocity bounded by $2\lambda\pi$.
We point out that because of the periodicity assumption on $f$,
and hence on $\QL(f,f)$, the collision operator $\QL(f,f)$ preserves in
time the mass contained in the period. This can be proved directly using the
property that $\B(-m,m)=0$. On the contrary, momentum and energy are not
preserved in time.

\begin{remark}
\rm
Note that, the Boltzmann equation itself is non conservative in a
bounded domain. Deterministic conservative schemes are usually constructed
by modifying the collision mechanics or the numerical method in order to
force the conservation properties in a bounded domain. However, it
is clear that even with this modification the accuracy of the solution
with respect to the original equation is guaranteed only if the size of
the velocity domain is large enough.
\end{remark}

Following the same strategy as in \cite{PR}, the spectral
accuracy of the method will be a consequence of a uniform
estimate on the true Boltzmann operator. This is the
content of the following proposition.

\begin{Prop} Let $B(q,\theta) = \Phi(|q|) b(\cos\theta)$, and
\[ Q_B(f,g) = \int_{\R^3} dv_* \int_{S^2} d\sigma\,
B(v-v_*,\theta) (g'_* f' - f_* f). \]
Then
\begin{equation}
\| Q_B(f,g) \|_{L^2} \leq C \|g\|_{L^1} \|f\|_{H^2},
\end{equation}
where $C$ is, up to a universal constant,
\[ \left [ \sup_{q \in \R^3} (1+|q|^2) \Phi(|q|) \right ]
\int_0^{\pi/2} d\theta\, \sin \theta b(\cos\theta) (1-\cos\theta). \]
\label{pr:stima}
\end{Prop}

\begin{remark}
\rm
\begin{itemize}
\item[]
\item[i)] By symmetrizing the kernel, we restricted the
domain of $\theta$ to $[0,\pi/2]$. This is always possible
thanks to the indiscernability of particles.
\item[ii)] In the applications, we shall use this proposition with
the support of $f$ and $g$ truncated over $[-\pi,\pi]^3$, so that the norm of
$g$ in $L^1$ is readily estimated by the norm in $L^2$, and
the constant $C$ is finite even for $\Phi(|q|) = |q|^\alpha$,
$\alpha >0$.
\item[iii)] In the {\em singular} case where $\Phi(|q|) = |q|^\alpha$,
$-3 < \gamma <0$, similar estimates hold, if one allows to replace
$\|g\|_{L^1}$ by (say) $\|g \ast \Phi \|_{L^\infty}$. This is
estimated by a norm of $g$ in some $L^p$ space, which in turn
can be estimated (at least locally) by the norm of $g$ in some
Sobolev space $H^\ell$.
\item[iv)] This estimate is clearly, in the case which is of interest for
us, uniform during the process of asymptotics of grazing
collisions.
\end{itemize}
\end{remark}

\begin{proof}
We use a duality argument~:
\[ \|Q_B(f,g)\|_{L^2} = \sup
\left \{ \int Q_B(f,g)\, \varphi, \qquad \|\varphi\|_{L^2} \leq 1 \right \}. \]
But
\[ \int Q_B(f,g)\varphi =
\int_{R^6} dv\, dv_*\, d\sigma\, B(v-v_*,\theta)
g_* f (\varphi' -\varphi) \]
\[ \leq \|g\|_{L^1}
\left [ \sup_{q \in \R^3} (1+|q|^2) \Phi(|q|) \right ]
\sup_{v_*\in \R^3}
\int_{R^3} \frac{dv}{1+|v-v_*|^2} \left |
\int_{S^2} d\sigma\, b(\cos\theta) f (\varphi' -\varphi) \right | \]
and
\begin{multline*}
 \int_{\R^3} dv\, \frac{f(v)}{1+|v-v_*|^2}
\left | \int d\sigma \, b(\cos\theta) (\varphi'-\varphi) \right | \\
\leq \|f\|_{H^2} \left \|\frac{1}{1+|v-v_*|^2}
\int_{S^2} d\sigma\, b(\cos\theta) (\varphi'-\varphi) \right \|_{H^{-2}}.
\end{multline*}

So it suffices to show that
\[ \left \| \frac{1}{1+|v-v_*|^2}
\int_{S^2} d\sigma\, b(\cos\theta) (\varphi'-\varphi) \right \|_{H^{-2}}
\leq C \|\varphi\|_{L^2}. \]
By translational invariance, it suffices to prove this estimate for
$v_*=0$. Since $1/(1+|v|^2)$ is a smooth function of $v$, it suffices to
prove that
\begin{equation} \label{goal}
\left \| \frac{1}{1+|v|^2} \int_{S^2} d\sigma\,
b(\theta) (\varphi'-\varphi) \right \|_{L^2}
\leq C \|\varphi\|_{H^2},
\end{equation}
where $\varphi'=\varphi((v+|v|\sigma)/2)$.

We now note that (up to a multiplicative constant)
\[ \int_{S^2} d\sigma\, b(\cos\theta) (\varphi'-\varphi)
= \int_0^{\pi/2} \sin\theta b(\cos\theta) \int_0^{2\pi} d\phi\,
(\varphi'-\varphi) \]
\begin{multline*}
= \int_0^{\pi/2} d\theta\, \sin\theta b(\cos\theta)
\nabla\varphi(v)\cdot (v'-v)\, d\phi \\
+ \int_0^{\pi/2} d\theta\, \sin\theta b(\cos\theta)
\int_0^{2\pi} \int_0^1 (1-t)
D^2\varphi (v+t(v'-v)) (v'-v,v'-v)\,dt\,d\phi.
\end{multline*}

Now, by symmetry,
\[ \int_0^{2\pi} d\phi\, \nabla \varphi(v)\cdot (v'-v)
= \nabla \varphi(v) \cdot v (\cos\theta -1), \]
and since $|v'-v|^2 = |v|^2 \sin^2(\theta/2)$, we find
\begin{multline*}
\left | \int_{S^2} d\sigma\, b(\cos\theta) (\varphi'-\varphi)
\right | \leq
C \left [ \int_0^{\pi/2} d\theta\, \sin\theta b(\cos\theta)
(1-\cos\theta) \right ] \\
\left [ |\nabla \varphi (v) | |v| + \sup_{\sigma\cdot v>0}
\int_0^1 dt\, (1-t) |D^2\varphi(v+t(v'-v))| |v|^2 \right ]
\end{multline*}
Then,~\eqref{goal} follows by noting that for each
$t \in [0,1]$, $\sigma \in S^2$ such that $\sigma\cdot v>0$,
one has $\int dv\, |D^2\varphi(v+t(v'-v))|^2
\leq C \int dv\, |D^2\varphi|^2\, dv$.
(the Jacobian of $v \mapsto v+t(v'-v)$ is bounded below).

\end{proof}

Next we state the consistency in the $L^2$-norm for the
approximation of the collision operator $\QL(f,f)$ with $\QL_N(f_N,f_N)$,
\begin{Thm}
Let $f \in H^2_p([-\pi,\pi]^3)$, then $\forall\, r\geq 0$
\begin{equation}
||\QL(f,f)-\QL_N(f_N,f_N)||_{2} \leq C \left(||f-f_N||_{H^2_p} +
\frac{||\QL(f_N,f_N)||_{H^r_p}}{N^r}\right),
\end{equation}
where $C$ depends on $||f||_2$.
\label{th:2}
\end{Thm}
\proof
First, we can split the error in two parts
\[
||\QL(f,f)-\QL_N(f_N,f_N)||_2 \leq ||\QL(f,f)-\QL(f_N,f_N)||_2
\]
\[
+ ||\QL(f_N,f_N)-\QL_N(f_N,f_N)||_2.
\]
Now clearly $\QL(f_N,f_N) \in \poly_{2N}$ and hence $\QL(f_N,f_N)$ is
periodic and infinitely smooth together with all its derivatives thus
\cite{CHQZ}
\begin{equation}
||\QL(f_N,f_N)-\QL_N(f_N,f_N)||_2 \leq \frac{C}{N^r} ||\QL(f_N,f_N)||_{H^r_p}, \quad
\forall\, r \geq 0.
\label{eq:Qb}
\end{equation}
By application of Proposition \ref{pr:stima} and of the identity
\[
\QL(f,f)-\QL(g,g) = \QL(f+g,f-g),
\]
we have
\[
||\QL(f,f)-\QL(f_N,f_N)||_2 = ||\QL(f+f_N,f-f_N)||_2 \leq C
||f+f_N||_1||f-f_N||_{H^2_p}
\]
\[
\leq 2 C_1 ||f||_2 ||f-f_N||_{H^2_p}.
\]
\endproof

Finally the following corollary states the uniform (in the grazing collision parameter) spectral accuracy of
the approximation of the collision operator
\begin{Cor}
Let $f \in H_p^r([-\pi,\pi]^3)$, $r\geq 2$ then
\begin{equation}
||\QL(f,f)-\QL_N(f_N,f_N)||_2 \leq \frac{C}{N^{r-2}} \left(||f||_{H^r_p} +
||\QL(f_N,f_N)||_{H^r_p}\right),
\end{equation}
\end{Cor}
\proof
It is enough to observe that
\[
||f-f_N||_{H^2_p} \leq \frac{C}{N^{r-2}} ||f||_{H^r_p}.
\]
\endproof
\begin{remark}
\rm
From the previous corollary it follows
\[
|<\QL(f,f),\varphi>-<\QL_N(f_N,f_N),\varphi>| \leq \frac{C}{N^{r-2}}||\varphi||_2 \left(||f||_{H^r_p} +
||\QL(f_N,f_N)||_{H^r_p}\right),
\]
and hence, by taking $\varphi=v,v^2$, the spectral accuracy of the moments.
\end{remark}

\section{Conclusions}

In this paper we have studied a Fourier spectral method that allows a numerical
passage from the Boltzmann equation to the FPL equation in the grazing
collision limit. In particular we have proved the uniform boundedness of
the kernel modes for any collision operator with a singular kernel,
provided the angular cross-section for momentum transfer (\ref{mt}) is bounded.
This permits to show that in the small grazing collision limit the method
provides a consistent discretization of the limiting FPL equation.
The uniform spectral accuracy of the method with respect to the grazing
collision parameter has been also given.
Moreover, we have derived an approximated formula (\ref{approx}) that
gives an intermediate asymptotic of the Boltzmann kernel modes which can be
evaluated with fast algorithms similarly to the FPL modes studied in~\cite{PRT}.
Numerical results based on this approximated formula and comparison with exact
solutions are actually under development and will be presented elsewhere.

\subsection*{Acknowledgements}
This research was partially supported by the grants ERBFMRXCT970157
(TMR-Network) from the EU, and by the National Reasearch Project (MURST), ``Numerical analysis: methods and scientific software''.
C.V. acknowledges the kind hospitality of the Departments of Mathematics of the Universities of Ferrara and Pavia, where the most part of the paper was done. G.T. acknowledges the support of the national Council for Researches, GNFM.

\newpage


\end{document}